\documentclass{svjour3}                     
\smartqed  
\usepackage{latexsym}
\usepackage{amssymb}
\usepackage{amsmath}
\usepackage{rotating}
\usepackage[misc]{ifsym}

\newtheorem{lem}{Lemma}[section]
\newtheorem{them}{Theorem}[section]

\newtheorem{pro}{Proposition}[section]

\journalname{Optim Letter}
\begin{document}

\title{On Metric Subregularity for Convex Constraint Systems by Primal Equivalent Conditions\thanks{This research was supported by the National Natural Science Foundation of P. R. China (grant  11401518), the Fok Ying-Tung Education Foundation (grant 151101) and the Scientific Research Foundation from Education Department of Yunnan Province under grant 2015Z009.}}

\titlerunning{Metric Subregularity by Primal Equivalent Conditions}



\author{Liyun Huang\and Zhou Wei}


\institute{Liyun Huang \at School of Mathematics and Information Science, Qujing Normal University,\\ Qujing 655011,Yunnan Province, People's Republic of China\\ \email{ynszlyb@163.com}
\and
 Zhou Wei(\Letter) \at Department of Mathematics, Yunnan University,
 Kunming 650091, People's Republic of China\\ \email{wzhou@ynu.edu.cn}}

\date{Received: date / Accepted: date}

\maketitle

\begin{abstract}
In this paper, we mainly study metric subregularity for a convex constraint system defined by a convex set-valued mapping and a convex constraint subset. The main work is to provide several primal equivalent conditions for metric subregularity by contingent cone and graphical derivative. Further it is proved that these primal equivalent conditions can characterize strong basic constraint qualification of convex constraint system given by Zheng and Ng [SIAM J. Optim., 18(2007), pp. 437-460].


\end{abstract}

\section{Introduction}
Many optimization problems appearing in variational analysis and mathematical programming can be modelled as finding a solution to a generalized equation. This generalized equation mathematically is defined as follows
\begin{equation}\label{1.1}
  \bar y\in F(x)
\end{equation}
where $F:X \rightrightarrows Y$ is a set-valued mapping between Banach spaces $X$ and $Y$ and $\bar y$ is a given point in $Y$.

It is well known that a key concept, when studying the behavior of the solution set to generalized equation \eqref{1.1}, is metric regularity. Recall that $F$ is said to be metrically regular at $\bar x\in F^{-1}(\bar y)$ if there exists $\tau\in (0,+\infty)$ such that
\begin{equation}\label{1.2}
  d(x,F^{-1}(y))\leq \tau d(y, F(x))\ \ {\rm for\ all}\ (x,y)\ {\rm close\ to} \ (\bar x,\bar y).
\end{equation}
The study of this concept can be traced back to the Robinson-Ursescu theorem, Lyusternik-Graves theorem or even Banach open mapping theorem. Readers are invited to consult \cite{A,BS,I,KK,Li,M1,M2} and references therein for many theoretical results on metric regularity and its various applications.

A weaker property of metric regularity is named as metric subregularity which means that inequality \eqref{1.2} holds only for fixed $\bar y$. Recall that $F$ is said to be metrically subregular at $\bar x\in F^{-1}(\bar y)$ if there exists $\tau\in (0,+\infty)$ such that
\begin{equation}\label{1.3}
  d(x,F^{-1}(\bar y))\leq \tau d(\bar y, F(x))\ \ {\rm for\ all}\ x\ {\rm close\ to} \ \bar x.
\end{equation}
This concept can be used to estimate the distance of a candidate $x$ to the solution set of generalized equation \eqref{1.1}. When the solution set $F^{-1}(\bar y)$ reduces to the singleton locally, \eqref{1.3} becomes to the strong metric subregularity of \eqref{1.1}; that is, $F$ is said to be strongly metrically subregular at $\bar x\in F^{-1}(\bar y)$ if there exists $\tau\in (0,+\infty)$ such that
$$
  \|x-\bar x\|\leq \tau d(\bar y, F(x))\ \ {\rm for\ all}\ x\ {\rm close\ to} \ \bar x.
$$
It is known that metric subregularity is closely related to calmness, error bounds, linear regularity and basic constraint qualification (BCQ) in optimization and has found a huge range of applications in areas of variational analysis and mathematical programming like  optimality conditions, variational inequalities, subdifferential theory, sensitivity analysis of generalized equations and convergence analysis of algorithms for solving equations or inclusions. For these reasons, the concept of metric subregularity has been extensively studied by many authors (see \cite{DR,Gf,HO,HW,IO,P1,W,WZ,Za,Z,Z1,Z2}).

In this paper, we mainly consider metric subregularity and targets at its primal criteria. In 2004, Dontchev and Rockafellar \cite{DR} studied metric subregularity of generalized equation \eqref{1.1} and proved that metric subregularity (at $\bar x$ for $\bar y\in F(\bar x)$) is equivalent to calmness of the inverse set-valued mapping $F^{-1}$ (at $(\bar y,\bar x)$). In 2007, Zheng and Ng \cite{Z1} discussed metric subregularity of a convex constraint system defined by a convex set-valued mapping and a convex constraint set. Using normal cone and coderivative, they introduced the concept of strong basic constraint qualification (strong BCQ) and established several dual characterizations in terms of strong BCQ for metric subregularity. Subsequently they \cite{Z2} considered metric subregularity for nonconvex generalized equation and gave its dual sufficient criteria.  Recently the authors \cite{HW} discussed several types of strong BCQs for nonconvex generalized equation and used these strong BCQs to provide necessary and/or sufficient conditions ensuring metric subregularity. To the best of our knowledge, these dual conditions for metric subregularity heavily rely on dual space properties, while one pretty natural idea on the study of metric subregularity is referring to primal space properties and sometimes it is convenient by this idea to deal with the case that there is little information on any dual space property like a metric space (see \cite{A,IO,K,K1,LM,MO,NT,P1}). Inspired by this observation, we mainly study primal criteria for metric subregularity of convex constraint system in this paper and aim to establish its several equivalent conditions by contingent cone and graphical derivative. To show clearly the primal spirit, we provide the self-contained proof for these equivalent conditions (see Section 3). We also prove that these primal equivalent conditions are necessary and sufficient for the strong BCQ given in \cite{Z1} (see Section 4).

The paper is organized as follows. In Section 2, we will give some preliminaries used in this paper. Our notation is basically standard and conventional in the area of variational analysis. Section 3 is devoted to main results on metric subregularity for a convex constraint system. Several primal equivalent conditions for  metric subregularity are established in terms of contingent cones and graphical derivatives. In Section 4, we mainly study the relation between primal equivalent conditions obtained in Section 3 and the strong BCQ of a convex constraint system. It is proved that these primal conditions can characterize the strong BCQ. The conclusion of this paper is presented in Section 5.

\section{Preliminaries}
Let $X, Y$ be Banach spaces with the closed unit balls denoted by
$B_{X}$ and $B_Y$. For a subset $C$ of $X$, let $\overline
C$ denote the closure of $C$. For any point $x\in X$ and $\delta>0$, let $B(x,\delta)$ denote an open ball with center $x$ and radius $\delta$.

Let $A$ be a closed convex subset of $X$ and $a\in A$. We denote by $T(A, a)$ the contingent (Bouligand) cone of $A$ at $a$ which is defined by
$$
T(A,a):=\mathop{\rm Limsup}\limits_{t\rightarrow 0^+}\frac{A-a}{t}.
$$
Thus, $v\in T(A,a)$ if and only if there exist $v_n\rightarrow v$ and $t_n\rightarrow 0^+$ such that $a+t_nv_n\in A$ for all $n$. It is easy to verify that $T(A, a)$ is a closed convex cone and
$$
T(A,a)=\overline{\bigcup_{t>0}\frac{A-a}{t}}.
$$

For a set-valued mapping $F:X\rightrightarrows Y$, we denote by
$${\rm gph}(F):=\{(x, y)\in X\times Y : y\in F(x)\}
$$
the graph of $F$. Recall that $F$ is said to be closed if ${\rm gph}(F)$ is a closed subset of $X\times Y$.

Let $F:X\rightrightarrows Y$ be a closed convex set-valued mapping and $(x,y)\in {\rm gph}(F)$. Recall that the graphical derivative $DF(x, y)$ of $F$ at $(x, y)$ is defined by
$$
DF(x, y)(u):=\{v\in Y: (u,v)\in T({\rm gph}(F),(x, y))\}.
$$
It is easy to verify that
$$
v\in DF(x, y)(u)\Longleftrightarrow u\in DF^{-1}(y,x)(v).
$$

We close this section with the following lemma cited from \cite[Lemma 2.1]{WZ}.

\begin{lem}
Let $X$ be a Banach space and  $\Omega$ be a nonempty closed convex
subset of $X$. Let $\gamma\in (0,\;1)$. Then for any $x\not\in \Omega$
there exists $z\in \Omega$ such that
$$\gamma\|x-z\|<d(x-z,T(\Omega,z)).$$
\end{lem}

\setcounter{equation}{0}

\section{Main results}

In this section, we mainly study metric subregularity of a convex constraint system and aim at providing its primal equivalent conditions in terms of contingent cone and graphical derivative.

Throughout the rest of this paper, we suppose that $F:X\rightrightarrows Y$ is a closed convex set-valued mapping and $A$  is a closed convex subset of $X$.

Let $\bar y\in Y$ be given. We consider the following convex constraint system:
\begin{equation}\label{3.1a}
  \bar y\in F(x)\ \ {\rm subject \ to} \ \ x\in A.
\end{equation}
We denote by $S:=F^{-1}(\bar y)\cap A$ the solution set of (3.1).

Recall from \cite{Z1} that convex constraint system (3.1) is said to be  metrically subregular at $\bar x\in S$ if there exist $\tau,\delta\in(0,+\infty)$ such that
\begin{equation}\label{3.1}
  d(x,S)\leq \tau(d(\bar y,F(x))+d(x, A))\ \  \forall x\in B(\bar x,\delta).
\end{equation}
We denote by ${\rm subreg}_AF(\bar x,\bar y)$ the modulus of metric subregularity at $\bar x$ for (3.1); that is,
\begin{equation}\label{3.3a}
  {\rm subreg}_AF(\bar x,\bar y):=\inf\{\tau>0: \exists \ \delta>0 \ {\rm s.t.} \ \eqref{3.1} \ {\rm holds}\},
\end{equation}
here we use the convention that the infimum over the empty set is $+\infty$. It is easy to verify that (3.1) is metrically subregular at $\bar x$ $\Longleftrightarrow$  ${\rm subreg}_AF(\bar x,\bar y)<+\infty$.

For convenience to provide primal equivalent conditions for metric subregularity, we consider two quantities which are defined by contingent cone and graphical derivative.

For any $x\in S$ and $\eta,\tau\in (0,+\infty)$ given, we study the following two relations:
\begin{equation}\label{3.2}
\begin{array}r
  DF^{-1}(\bar y,x)(\eta_1B_Y)\cap(T(A,x)+\eta_2B_X)\subset T(S,x)+B_X,  \forall \eta_1,\eta_2\geq 0\  {\rm with} \\  \eta_1+\eta_2<\eta,
\end{array}
\end{equation}
and
\begin{equation}\label{3.13a}
   d(h, T(S, x))\leq \tau (d(0, DF(x,\bar y)(h))+d(h, T(A, x))) \ \ \forall h\in X.
\end{equation}

Let $\bar x\in S$. In terms of these two relations, we define the following quantities:
\begin{equation}\label{3-06}
 \eta_A(F, \bar x,\bar y):=\sup\{\eta>0: \exists \ \delta>0 \ {\rm s.t.} \  \eqref{3.2}\ {\rm holds\ for \ any}\ x\in S\cap B(\bar x,\delta)\},
\end{equation}
and
\begin{equation}\label{3-07}
  \tau_A(F, \bar x,\bar y):=\inf\{\tau>0: \exists \ \delta>0 \ {\rm s.t.} \ \eqref{3.13a}\ {\rm holds\ for \ any}\  x\in S\cap B(\bar x,\delta)\},
\end{equation}
where the supremum over the empty set is $0$.

The following theorem, as one main result in the paper, is to establish an accurate quantitative relation among modulus \eqref{3.3a}, quantities \eqref{3-06} and \eqref{3-07}. This theorem also provides primal characterizations for metric subregularity of convex constraint system (3.1).
\begin{them}
Let $\bar x\in S$. Then
\begin{equation}\label{3-08}
\frac{1}{{\rm subreg}_AF(\bar x, \bar y)}=\eta_A(F, \bar x,\bar y)=\frac{1}{\tau_A(F,\bar x, \bar y)}.
\end{equation}

\end{them}

\noindent{\it Proof} We first consider the case that ${\rm subreg}_AF(\bar x, \bar y)<+\infty$. We prove that
\begin{equation}\label{3-9}
  \frac{1}{{\rm subreg}_AF(\bar x, \bar y)}=\eta_A(F, \bar x,\bar y).
\end{equation}

Let $\tau>{\rm subreg}_AF(\bar x, \bar y)$. Then there exists $\delta>0$ such that \eqref{3.1} holds. Let $\eta\in (0,\frac{1}{\tau})$ and $x\in S\cap B(\bar x,\delta)$. Take any $\eta_1,\eta_2\in [0, +\infty)$ such that $\eta_1+\eta_2<\eta$ and any $u\in  DF^{-1}(\bar y,x)(\eta_1B_Y)\cap(T(A,x)+\eta_2B_X)$. Then there exist $z\in B_Y,b\in B_X$, $t_n\rightarrow 0^+,s_n\rightarrow 0^+$, $(v_n,u_n)\rightarrow (\eta_1z,u)$ and $w_n\rightarrow u+\eta_2b$ such that
\begin{equation}\label{3.4}
  x+t_nu_n\in F^{-1}(\bar y+t_nv_n)\ \ {\rm and} \ \ x+s_nw_n\in A\ \ \forall n\in \mathbb{N}.
\end{equation}
Since $F$ is a convex set-valued mapping and $A$ is convex, it follows from \eqref{3.4} that there exist $\lambda_n\rightarrow 0^+$ such that
\begin{equation*}
  x+\lambda_nu_n\in F^{-1}(\bar y+\lambda_nv_n)\ \ {\rm and} \ \ x+\lambda_nw_n\in A\ \ \forall n\in \mathbb{N}.
\end{equation*}
This and \eqref{3.1} imply that for any $n$ sufficiently large, one has
\begin{eqnarray*}
d(x+\lambda_nu_n, S)&\leq&\tau(d(\bar y, F(x+\lambda_nu_n))+d(x+\lambda_nu_n, A))\\
&\leq&\tau\lambda_n(\|v_n\|+\|u_n-w_n\|),
\end{eqnarray*}
and consequently
$$
d(u_n,T(S,x))\leq d\Big(u_n,\frac{S-x}{\lambda_n}\Big)\leq\tau
(\|v_n\|+\|u_n-w_n\|)
$$
(thanks to the convexity of $F$ and $A$). Taking limits as $n\rightarrow\infty$, one has
$$
d(u,T(S,x))\leq\tau
(\|\eta_1z\|+\|u-(u+\eta_2b)\|)\leq \tau(\eta_1+\eta_2).
$$
Thus,
$$
u\in\overline{T(S,x)+\tau(\eta_1+\eta_2)B_X}\subset T(S,x)+B_X
$$
thanks to $\eta_1+\eta_2<\eta<\frac{1}{\tau}$. This means that $\eta_A(F,\bar x,\bar y)\geq \eta$. By taking limits as $\eta \uparrow \frac{1}{\tau}$ and then $\tau\downarrow{\rm subreg}_AF(\bar x, \bar y)$, one has
\begin{equation}\label{3-11}
  \eta_A(F, \bar x,\bar y)\geq\frac{1}{{\rm subreg}_AF(\bar x, \bar y)}>0.
\end{equation}

Let $\eta\in (0, \eta_A(F, \bar x,\bar y))$. Then there exists $\delta>0$ such that \eqref{3.2} holds for all $x\in S\cap B(\bar x,\delta)$. Let $\tau>\frac{1}{\eta}$. We claim that
\begin{equation}\label{3.6}
  d(x, S)\leq\tau(d(\bar y, F(x))+d(x, A))\ \ \forall x\in B(\bar x, \frac{\delta}{2}).
\end{equation}
Granting this, it follows that
$$
\frac{1}{{\rm subreg}_AF(\bar x, \bar y)}\geq\frac{1}{\tau}.
$$
By taking limits as $\tau\downarrow\frac{1}{\eta}$ and then $\eta\uparrow \eta_A(F, \bar x,\bar y)$, one has
$$
\frac{1}{{\rm subreg}_AF(\bar x, \bar y)}\geq\eta_A(F, \bar x,\bar y).
$$
This and \eqref{3-11} imply that \eqref{3-9} holds.

It is clear that \eqref{3.6} holds for any $x\in B(\bar x, \frac{\delta}{2})\cap S$. Let $x\in B(\bar x, \frac{\delta}{2})\backslash S$. Then $d(x, S)\leq\|x-\bar x\|<\frac{\delta}{2}$. Take any $\gamma\in (\frac{2d(x, S)}{\delta}, 1)$. By virtue of Lemma 2.1, there exists $z\in S$ such that
\begin{equation}\label{3.7}
  \gamma\|x-z\|\leq d(x-z, T(S,z)).
\end{equation}
Since $S-z\subset T(S,z)$ by the convexity of $S$, it follows that
\begin{equation*}
  \gamma\|x-z\|\leq d(x-z, T(S,z))\leq d(x, S).
\end{equation*}
From this and the choice of $\gamma$, one has $\|x-z\|\leq \frac{d(x, S)}{\gamma}<\frac{\delta}{2}$ and thus
\begin{equation}\label{3.8a}
  \|z-\bar x\|\leq \|z-x\|+\|x-\bar x\|<\frac{\delta}{2}+\frac{\delta}{2}=\delta.
\end{equation}
Choose sequences $\{y_n\}\subset F(x)$ and $\{u_n\}\subset A$ such that
\begin{equation}\label{3.8}
  \|\bar y-y_n\|\rightarrow d(\bar y, F(x))\ \ {\rm and} \ \ \|x-u_n\|\rightarrow d(x, A).
\end{equation}
For each $n\in\mathbb{N}$, we denote
$$
r_n:=\frac{1}{\tau(\|\bar y-y_n\|+\|x-u_n\|)}>0
$$
and take $\varepsilon_n>0$ sufficiently small such that
\begin{equation}\label{3.9}
  r_n(\|\bar y-y_n\|+(1+\varepsilon_n)\|x-u_n\|)<\eta.
\end{equation}
Note that $\|r_n(x-z)-r_n(u_n-z)\|=r_n\|x-u_n\|$ and so
$$
d(r_n(x-z),T(A,z))\leq\|r_n(x-z)-r_n(u_n-z)\|=r_n\|x-u_n\|.
$$
This implies that
\begin{equation}\label{3.10}
  r_n(x-z)\in\overline{T(A,z)+r_n\|x-u_n\|B_X}\subset T(A,z)+(1+\varepsilon_n)r_n\|x-u_n\|B_X.
\end{equation}
Since
$$
r_n(x-z)\in DF^{-1}(\bar y,z)(r_n(y_n-\bar y))\subset DF^{-1}(\bar y,z)(r_n\|y_n-\bar y\|B_Y),
$$
it follows from \eqref{3.2}, \eqref{3.8a}, \eqref{3.9} and \eqref{3.10} that
$$
r_n(x-z)\in T(S,z)+B_X
$$
and consequently
$$
d(x-z, T(S,z))\leq\frac{1}{r_n}=\tau(\|\bar y-y_n\|+\|x-u_n\|).
$$
Taking limits as $n\rightarrow \infty$, one has
$$
d(x-z, T(S,z))\leq \tau(d(\bar y, F(x))+d(x, A)).
$$
This and \eqref{3.7} mean that
$$
\gamma d(x, S)\leq \tau(d(\bar y, F(x))+d(x, A)).
$$
Taking limits as $\gamma\rightarrow 1^-$, one has
$$
d(x, S)\leq \tau(d(\bar y, F(x))+d(x, A)).
$$
Hence \eqref{3.6} holds.

We next prove that
\begin{equation}\label{3-19}
\eta_A(F, \bar x,\bar y)=\frac{1}{\tau_A(F, \bar x, \bar y)}.
\end{equation}

Let $\eta\in (0, \eta_A(F, \bar x,\bar y))$. Then there exists $\delta>0$ such that \eqref{3.2} holds for all $x\in S\cap B(\bar x,\delta)$. Let $\tau>\frac{1}{\eta}$ and $x\in S\cap B(\bar x,\delta)$. Take any $h\in X$. Choose sequences $\{v_n\}\subset DF(x,\bar y)(h)$ and $\{u_n\}\subset T(A,x)$ such that
\begin{equation*}
  \|v_n\|\rightarrow d(0, DF(x,\bar y)(h))\ \ {\rm and }\ \ \|h-u_n\|\rightarrow d(h,T(A,x)).
\end{equation*}
For each $n\in\mathbb{N}$, let
$$
r_n:=\frac{1}{\tau(\|v_n\|+\|h-u_n\|)}>0.
$$
Noting that $h=u_n+\|h-u_n\|\cdot\frac{h-u_n}{\|h-u_n\|}$ and $v_n\in DF(x,\bar y)(h)$, it follows that
\begin{eqnarray*}
r_nh&\in& DF^{-1}(\bar y,x)(r_nv_n)\cap\big( T(A, x)+r_n\|h-u_n\|\cdot\frac{h-u_n}{\|h-u_n\|}\big)\\
&\subset&DF^{-1}(\bar y,x)(r_n\|v_n\|B_Y)\cap( T(A, x)+r_n \|h-u_n\|B_X).
\end{eqnarray*}
Since $r_n(\|v_n\|+\|h-u_n\|)=\frac{1}{\tau}<\eta$, by \eqref{3.2}, one has
$$
r_nh\in T(S,x)+B_X
$$
and consequently
$$
d(r_nh, T(S,x))\leq 1.
$$
This implies that
$$
d(h, T(S,x))\leq \frac{1}{r_n}=\tau(\|v_n\|+\|h-u_n\|).
$$
Taking limits as $n\rightarrow \infty$, one has that $\tau_A(F, \bar x, \bar y)\leq \tau$ and thus
\begin{equation*}\label{3-20}
  \frac{1}{\tau_A(F, \bar x, \bar y)}\geq\eta_A(F, \bar x,\bar y)
\end{equation*}
by taking limits as $\tau\downarrow\frac{1}{\eta}$ and then $\eta \uparrow \eta_A(F, \bar x, \bar y)$.

Let $\tau>\tau_A(F,\bar x,\bar y)$. Then there exists $\delta>0$ such that \eqref{3.13a} holds for any $x\in S\cap B(\bar x,\delta)$. Let $\eta\in (0,\frac{1}{\tau})$ and $x\in S\cap B(\bar x,\delta)$. Take any $\eta_1,\eta_2\in[0,+\infty)$ with $\eta_1+\eta_2<\eta$ and $u\in DF^{-1}(\bar y,x)(\eta_1B_Y)\cap(T(A, x)+\eta_2B_X)$. Then there exist $(v,b)\in B_Y\times B_X$ such that
$$
u\in DF^{-1}(\bar y,x)(\eta_1v)\ \ {\rm and} \ \ u-\eta_2b\in T(A,x).
$$
By virtue of \eqref{3.13a}, one has
\begin{eqnarray*}
d(u,T(S,x))&\leq& \tau (d(0, DF(x,\bar y)(u))+d(u, T(A, x)))\\
&\leq&\tau(\eta_1\|v\|+\eta_2\|b\|)\\
&\leq& \tau(\eta_1+\eta_2).
\end{eqnarray*}
Thus
$$
u\in \overline{T(S,x)+\tau(\eta_1+\eta_2)B_X}\subset T(S,x)+B_X
$$
as $\eta_1+\eta_2<\eta<\frac{1}{\tau}$. This means that $\eta_A(F,\bar x,\bar y)\geq \eta$ and consequently
$$
  \eta_A(F, \bar x,\bar y)\geq\frac{1}{\tau_A(F, \bar x, \bar y)}
$$
by taking limits as $\eta\uparrow\frac{1}{\tau}$ and then $\tau \downarrow \tau_A(F, \bar x, \bar y)$. Hence \eqref{3-19} holds.

Next, we consider the case that ${\rm subreg}_AF(\bar x, \bar y)=+\infty$. Then $\eta_A(F,\bar x,\bar y)=0$ (Otherwise, one can verify that ${\rm subreg}_AF(\bar x, \bar y)\leq\frac{1}{\eta_A(F,\bar x,\bar y)}<+\infty$, which is a contradiction).  From this, we can further prove that $\tau_A(F, \bar x, \bar y)=+\infty$. Hence \eqref{3-08} holds. The proof is complete.\hfill$\Box$\\

\noindent{\bf Remark 3.1.}  Note that the concept of metric $q$-subregularity ($q>0$) for generalized equations has been well studied (cf. \cite{K,K1,LM,MO} and references therein). Constraint system \eqref{3.1a} reduces to generalized equation $\bar y\in F(x)$ if $A=X$. Li and Mordukhovich \cite{LM} used modulus estimate via coderivatives for ensuring metric $q$-subregularity of $\bar y\in F(x)$ (see \cite[Theorem 3.3]{LM}). Mordukhovich and Ouyang \cite{MO}  studied metric $q$-subregularity of generalized equation $\bar y\in F(x)$ in which $Y=X^*$ and $F=\partial f$ (the subdifferential of  function $f$ on $X$). It is proved that metric $q$-subregularity of subdifferential $\partial f$ can be characterized by higher order growth conditions of $f$ (see \cite[Theorem 3.4]{MO}).\hfill$\Box$\\

For the special case, we have sharper primal results on metric subregularity; that is, the validity of \eqref{3.2} only at $\bar x$ can ensure the metric subregularity of convex constraint system (3.1).

\begin{them}
Let $\bar x\in S$. Suppose that there exists a closed convex cone $K$ and a neighborhood $V$ of $\bar x$ such that $S\cap V=(\bar x+K)\cap V$.  Then
\begin{equation}\label{3.20c}
\frac{1}{{\rm subreg}_AF(\bar x, \bar y)}=\sup\{\eta>0: \eqref{3.2}\ {\it holds\ only \ with}\ x=\bar x\}.
\end{equation}
\end{them}

\noindent{\it Proof} We denote
$$
\alpha:=\sup\{\eta>0: \eqref{3.2}\ {\rm holds\ only \ with}\ x=\bar x\}.
$$
By Theorem 3.1, one has
\begin{equation}
\frac{1}{{\rm subreg}_AF(\bar x, \bar y)}=\eta_A(F,\bar x, \bar y)\leq \alpha.
\end{equation}
Let $\eta\in (0, \alpha)$ and $\tau>\frac{1}{\eta}$. Take $\delta>0$ such that $B(\bar x,2\delta)\subset V$. Then
$$
S\cap B(\bar x,2\delta)=(\bar x+K)\cap B(\bar x,2\delta).
$$
Since $K$ is a convex cone, one can verify that
\begin{equation}\label{3.14}
  T(S,\bar x)=K.
\end{equation}
Let $u\in B(\bar x,\delta)\backslash S$. Take any $\gamma\in(\frac{d(u,S)}{\delta}, 1)$. By Lemma 2.1, there exists  $z\in S$ such that
\begin{equation}\label{3.15}
  \gamma\|u-z\|\leq d(u-z,T(S,z)).
\end{equation}
Note that $S-z\subset T(S,z)$ by the convexity of $S$ and so $\|u-z\|\leq \frac{d(u-z,T(S,z))}{\gamma}<\delta$. Then
$$
\|z-\bar x\|\leq \|z-u\|+\|u-\bar x\|<\delta+\delta=2\delta.
$$
Noting that $K$ is a convex cone, it follows from \eqref{3.14} that
\begin{equation*}
  T(S,z)=T(\bar x+K,z)=T(K,z-\bar x)\supset  K-z+\bar x=T(S,\bar x)-z+\bar x.
\end{equation*}
This and \eqref{3.15} imply that
\begin{equation}\label{3.16}
  \gamma d(u, S)\leq \gamma\|u-z\|\leq d(u-\bar x, T(S,\bar x)).
\end{equation}
Similar to the proof of $
\frac{1}{{\rm subreg}_AF(\bar x, \bar y)}\geq\eta_A(F, \bar x,\bar y)$ in Theorem 3.1, by using \eqref{3.16}, one can prove that
$$
d(u, S)\leq\tau(d(\bar y, F(u))+d(u, A)).
$$
This means that ${\rm subreg}_AF(\bar x, \bar y)\leq\tau$ and thus by taking limits as $\tau\downarrow\frac{1}{\eta}$ and then $\eta \uparrow \alpha$, one has that
$$
\frac{1}{{\rm subreg}_AF(\bar x, \bar y)}\geq\alpha.
$$
Hence \eqref{3.20c} holds. The proof is complete. \hfill$\Box$\\

Next, we consider strong metric subregularity of convex constraint system. Recall that convex constraint system (3.1) is said to be strongly metrically subregular at $\bar x\in S$ if there exist $\tau,\delta\in(0,+\infty)$ such that
\begin{equation}\label{3.23c}
\|x-\bar x\|\leq\tau(d(\bar y, F(x))+d(x,A)) \ \ \forall x\in B(\bar x,\delta).
\end{equation}
We denote by ${\rm ssubreg}_AF(\bar x,\bar y)$ the modulus of strong metric subregularity at $\bar x$ for (3.1); that is,
\begin{equation}\label{3.24c}
  {\rm ssubreg}_AF(\bar x,\bar y):=\inf\{\tau>0: \exists \delta>0\ {\rm s.t.} \ \eqref{3.23c} \ {\rm holds} \}.
\end{equation}
Thus, (3.1) is strongly metrically subregular at $\bar x$ $\Longleftrightarrow$  ${\rm ssubreg}_AF(\bar x,\bar y)<+\infty$.

It is easy to verify that (3.1) is strongly metrically subregular at $\bar x\in S$ if and only if (3.1) is metrically subregular at $\bar x$ and $S=\{\bar x\}$.

To present primal characterizations for strong metric subregularity of (3.1), we study one inclusion which is given by contingent cone and graphical derivative.

For $\bar x\in S$ and $\eta\in (0,+\infty)$ given, we consider the following inclusion:
\begin{equation}\label{3.19a}
  DF^{-1}(\bar y,\bar x)(\eta_1B_Y)\cap(T(A,\bar x)+\eta_2B_X)\subset B_X\ \ \forall \eta_1,\eta_2\geq 0\ {\rm with} \ \eta_1+\eta_2<\eta.
\end{equation}

\begin{pro}
Let $\bar x\in S$ and $\eta\in (0,+\infty)$ be such that \eqref{3.19a} holds. Then $S=\{\bar x\}$.
\end{pro}
\noindent{\it Proof} Let $x\in S$. Then for any $t>0$, one has
$$
t(x-\bar x)\in DF^{-1}(\bar y,\bar x)(0)\cap T(A,\bar x)
$$
as ${\rm gph}(F)$ and $A$ are convex. It follows from \eqref{3.19a} that
$$
t(x-\bar x)\in B_X\ \ \forall t>0.
$$
This implies that $x=\bar x$ and thus $S=\{\bar x\}$ holds. The proof is complete. \hfill$\Box$\\

The following theorem provides an accurate quantitative estimate on the modulus of \eqref{3.24c} and also gives primal characterizations for strong metric subregularity of (3.1). The proof can be obtained by Theorem 3.1 and Proposition 3.1.

\begin{them}
Let $\bar x\in S$. Then
\begin{equation}\label{3.26}
\frac{1}{{\rm ssubreg}_AF(\bar x, \bar y)}=\sup\{\eta>0: \eqref{3.19a}\ {\it holds} \}.
\end{equation}
\end{them}

\begin{pro}
Let $\bar x\in S$. Suppose that $DF^{-1}(\bar y,\bar x)(B_Y)\cap(T(A,\bar x)+B_X)\cap B_X$ is relatively compact. Then ${\rm ssubreg}_AF(\bar x,\bar y)<+\infty$ if and only if
\begin{equation}\label{3.22}
  DF^{-1}(\bar y,\bar x)(0)\cap T(A,\bar x)=\{0\}.
\end{equation}

\end{pro}

\noindent{\it Proof} The necessity part. By Theorem 3.3, there exists $\eta>0$ such that \eqref{3.19a} holds. Then
$$
 DF^{-1}(\bar y,\bar x)(0)\cap T(A,\bar x)\subset \varepsilon B_{X} \ \ \forall \varepsilon>0.
$$
This means that \eqref{3.22} holds.

The sufficiency part. By using Theorem 3.3, we only need to prove that there exists $\eta>0$ such that \eqref{3.19a} holds.

Suppose on the contrary that for any $n\in \mathbb{N}$, there exist $r_n,s_n\geq 0$ with $r_n+s_n<\frac{1}{n}$, $(y_n,b_n)\in B_Y\times B_X$ and $x_n\in X$ such that
\begin{equation}\label{3.25}
  x_n\in\big(DF^{-1}(\bar y,\bar x)(r_ny_n)\cap (T(A,\bar x)+s_nb_n)\big)\backslash B_X.
\end{equation}
This implies that
\begin{eqnarray*}
\frac{x_n}{\|x_n\|}\in DF^{-1}(\bar y,\bar x)(B_Y)\cap (T(A,\bar x)+B_X).
\end{eqnarray*}
Since $DF^{-1}F(\bar y,\bar x)(B_Y)\cap(T(A,\bar x)+B_X)\cap B_X$ is relatively compact, without loss of generality, we can assume that $\frac{x_n}{\|x_n\|}\rightarrow x_0$ with $\|x_0\|=1$ (considering subsequence if necessary). By \eqref{3.25}, one has
$$
\frac{x_n}{\|x_n\|}\in DF^{-1}(\bar y,\bar x)(\frac{r_n}{\|x_n\|}y_n)\cap (T(A,\bar x)+\frac{s_n}{\|x_n\|}b_n).
$$
By taking the limit as $n\rightarrow \infty$, one has
$$
x_0\in  DF^{-1}(\bar y,\bar x)(0)\cap T(A,\bar x),
$$
which contradicts \eqref{3.22} since $x_0\not=0$. The proof is complete. \hfill$\Box$\\

\setcounter{equation}{0}

\section{Equivalence with strong BCQ of convex constraint systems}

In this section, we mainly study interrelationship between primal equivalent conditions in Section 3 and strong BCQ of convex constraint system (3.1) introduced by Zheng and Ng \cite{Z1}.  The main result in this section shows that primal equivalent conditions obtained can characterize the strong BCQ of convex constraint system (3.1). We recall some definitions and notations.

Let $X^*, Y^*$ denote the dual spaces of $X$
and $Y$ respectively. For a closed convex subset $A$ of $X$ and a point $a\in A$, we denote by $\delta_A$ the indicator function of $A$ and denote by $N(A,a)$ the normal cone of $A$ at $a$ which is defined as
$$
N(A,a):=\{x^*\in X^*:\;\langle x^*,x-a\rangle\leq0\;\;\forall x\in
A\}.$$
For a closed convex set-valued mapping $F:X\rightrightarrows Y$, let $D^*F(x, y): Y^*\rightrightarrows X^*$ denote the coderivative of $F$ at $(x, y)\in{\rm gph}(F)$ which is defined as
$$
D^*F(x, y)(y^*):=\{x^*\in X^* : (x^*, -y^*)\in N({\rm gph}(F), (x, y))\}\ \ \forall y^*\in Y^*.
$$
For a proper lower semicontinuous convex function $\psi: X\rightarrow \mathbb{R}\cup\{+\infty\}$, recall that the subdifferential of $\psi$ at $\bar x\in{\rm dom}(\psi):=\{x\in X: \psi(x)<+\infty\}$, denoted by $\partial\psi(\bar x)$, is defined by
$$
\partial\psi(\bar x):=\{x^*\in X^*: \langle x^*,x-\bar x\rangle\leq\psi(x)-\psi(\bar x)\ \ \forall x\in X\}.
$$

Recall from \cite{Z1} that convex constraint system (3.1) is said to have the strong BCQ at $x\in S$ if there exists $\tau\in(0,+\infty)$ such that
\begin{equation}\label{4.1}
  N(S,x)\cap B_{X^*}\subset\tau(D^*F(x,\bar y)(B_{Y^*})+N(A,x)\cap B_{X^*}).
\end{equation}

It is known that Zheng and Ng \cite{Z1} proved dual characterizations for metric subregularity of a convex constraint system by strong BCQ; that is, convex constraint system (3.1) is metrically subregular at $\bar x\in S$ if and only if there exists $\tau>0$ such that \eqref{4.1} holds for all $x\in S$ close to $\bar x$ with the same constant $\tau$.\\

We are now in a position to investigate the equivalent relation between primal space property of \eqref{3.2} and strong BCQ of \eqref{4.1} for convex constraint system (3.1). The following proposition shows the accurate quantitative relation between them.

\begin{pro}
Let $\bar x\in S$.  Then
\begin{equation}\label{4-2}
  \frac{1}{\eta_A(F, \bar x,\bar y)}=\inf\{\tau>0: \exists \ \delta>0 \ {\it s.t.} \ \eqref{4.1}\ {\it holds\ for \ any}\  x\in S\cap B(\bar x,\delta)\}.
\end{equation}

\end{pro}

\noindent{\it Proof} We denote that
$$
\beta:=\inf\{\tau>0: \exists \ \delta>0 \ {\rm s.t.} \ \eqref{4.1}\ {\rm holds\ for \ any}\  x\in S\cap B(\bar x,\delta)\}.
$$
Let $\tau>\beta$. Then there exists $\delta>0$ such that \eqref{4.1} holds for all $x\in S\cap B(\bar x,\delta)$. Let $\eta\in (0,\frac{1}{\tau})$ and $x\in S\cap B(\bar x,\delta)$. Take any $\eta_1,\eta_2\geq 0$ with $\eta_1+\eta_2<\eta$. Choose $\varepsilon>0$ sufficiently small such that
\begin{equation}\label{4.3}
  (1+\varepsilon)(\eta_1+\eta_2)<\eta.
\end{equation}
We claim that
\begin{equation}\label{4.4}
  DF^{-1}(\bar y,x)((1+\varepsilon)\eta_1B_Y)\cap(T(A,x)+(1+\varepsilon)\eta_2B_X)\subset\overline{ T(S,x)+B_X}.
\end{equation}
Granting this, it follows that
$$
  DF^{-1}(\bar y,x)((1+\varepsilon)\eta_1B_Y)\cap(T(A,x)+(1+\varepsilon)\eta_2B_X)\subset T(S,x)+(1+\varepsilon)B_X
$$
and thus \eqref{3.2} holds. This means that $\eta_A(F,\bar x,\bar y)\geq\eta$ and
\begin{equation}\label{4-5}
\frac{1}{\eta_A(F,\bar x,\bar y)}\leq\beta
\end{equation}
by taking limits as $\eta\uparrow \frac{1}{\tau}$ and then $\tau\downarrow\beta$.

Suppose on the contrary that \eqref{4.4} does not hold. Then there exists $x_0\in X$ such that
\begin{equation}\label{4.5}
  x_0\in \big( DF^{-1}(\bar y,x)((1+\varepsilon)\eta_1B_Y)\cap(T(A,x)+(1+\varepsilon)\eta_2B_X)\big)\backslash\overline{ T(S,x)+B_X}.
\end{equation}
By the separation theorem, there exists $x^*_0\in X^*$ with $\|x_0^*\|=1$ such that
\begin{equation}\label{4.6}
  \langle x^*_0,x_0\rangle>\sup\{\langle x_0^*,u \rangle: u\in T(S,x)+B_X\}=1.
\end{equation}
This implies that $x^*_0\in N(S,x)\cap B_{X^*}$. By virtue of the strong BCQ of \eqref{4.1}, there exist $x_1^*\in D^*F(x,\bar y)(y_1^*)$ for some $y^*_1\in B_{Y^*}$ and $x_2^*\in N(A,x)\cap B_{X^*}$ such that
\begin{equation}\label{4.7}
  x_0^*=\tau(x_1^*+x_2^*).
\end{equation}
By using \eqref{4.5}, there exist $v_0\in (1+\varepsilon)\eta_1B_Y$ and $u_0\in(1+\varepsilon)\eta_2B_X$ such that
$$
x_0\in  DF^{-1}(\bar y,x)(v_0)\ \ {\rm and}  \ \ x_0-u_0\in T(A,x).
$$
Then,
$$
\langle x_1^*,x_0\rangle\leq \langle y_1^*, v_0\rangle\leq (1+\varepsilon)\eta_1
$$
and
$$
 \langle x_2^*,x_0\rangle=\langle x_2^*,x_0-u_0+u_0\rangle\leq\langle x_2^*,u_0\rangle\leq (1+\varepsilon)\eta_2.
$$
This and \eqref{4.3} imply that
$$
\langle x_0^*,x_0\rangle=\tau(\langle x_1^*,x_0\rangle+\langle x_2^*,x_0\rangle)\leq \tau(1+\varepsilon)(\eta_1+\eta_2)<1,
$$
which contradicts \eqref{4.6}. Hence \eqref{4.4} holds.

Let $\eta\in (0, \eta_A(F, \bar x,\bar y))$. Then there exists $\delta>0$ such that \eqref{3.2} holds for all $x\in S\cap B(\bar x,\delta)$. Let $\tau>\frac{1}{\eta}$ and $x\in S\cap B(\bar x,\delta)$. Take any $x^*\in N(S,x)\cap B_{X^*}$. Note that
$$
N(S,x)\cap B_{X^*}=N(T(S,x),0)\cap B_{X^*}=\partial d(\cdot,T(S,x))(0).
$$
Then for any $h\in X$ and any $v\in DF(x,\bar y)(h)$, by Theorem 3.1, one has
\begin{eqnarray*}
\langle x^*, h\rangle\leq d(h,T(S,x))&\leq&\tau (d(0, DF(x,\bar y)(h))+d(h, T(A, x)))\\
&\leq&\tau(\|v\|+d(h, T(A, x))).
\end{eqnarray*}
Let $\Psi:X\times Y\rightarrow \mathbb{R}\cup\{+\infty\}$ be defined as
$$
\Psi(h,v):=\delta_{T({\rm gph}(F),(x,\bar y))}(h,v)+d(h, T(A, x))+\|v\|\ \ \forall (h,v)\in X\times Y.
$$
Then
\begin{equation}\label{4.8a}
  \langle(\frac{x^*}{\tau},0), (h,v)\rangle\leq\Psi(h,v)-\Psi(0,0) \ \ \forall (h,v)\in X\times Y.
\end{equation}
Noting that $\Psi$ is a proper lower semicontinuous convex function, it follows from \eqref{4.8a} that
\begin{eqnarray*}
(\frac{x^*}{\tau},0)\in\partial \Psi(0,0)&=&N(T({\rm gph}(F),(x,\bar y)), (0,0))+\partial d(\cdot,T(A,x))(0)\times B_{Y^*}\\
&=&N({\rm gph}(F),(x,\bar y))+(N(A,x)\cap B_{X^*})\times B_{Y^*}
\end{eqnarray*}
(the first equation follows from \cite[Theorem 3.16]{Ph}). This implies that
$$
x^*\in\tau(D^*F(x,\bar y)(B_{Y^*})+N(A,x)\cap B_{X^*})
$$
and consequently
$$
 N(S,x)\cap B_{X^*}\subset\tau(D^*F(x,\bar y)(B_{Y^*})+N(A,x)\cap B_{X^*}).
$$
This means that \eqref{4.1} holds and thus $\beta\leq\tau$. By taking limits as $\tau\downarrow\frac{1}{\eta}$ and then $\eta\uparrow\eta_A(F,\bar x,\bar y)$, one has
$$
\beta\leq\frac{1}{\eta_A(F,\bar x,\bar y)}.
$$
This and \eqref{4-5} imply that \eqref{4-2} holds. The proof is complete.\hfill$\Box$\\

\noindent{\bf Remark 4.1.} Given $\bar x\in S$,  by Theorem 3.1 and Proposition 4.1, one has that
$$
  \tau_A(F, \bar x,\bar y)=\inf\{\tau>0: \exists \ \delta>0 \ {\rm s.t.} \ \eqref{4.1}\ {\rm holds\ for \ any}\  x\in S\cap B(\bar x,\delta)\}.
$$
This shows that primal space property \eqref{3.13a} can also be used to characterize the strong BCQ of \eqref{4.1}.

\section{Conclusions}
This paper is devoted to the study of metric subregularity for a convex constraint system. Compared with the literature in dealing with metric subregularity by using dual tools like normal cone, subdifferential or coderivative, several primal equivalent conditions are established in terms of contingent cones and graphical derivatives. It is proved that these primal equivalent conditions can characterize the strong BCQ of the convex constraint system given in \cite{Z1}. This also demonstrates that the strong BCQ is essentially the dual space counterpart of these primal space properties. \\

\noindent{\bf Acknowledgment.} The authors wish to thank two anonymous referees for their many valuable comments and suggestions which help us to improve the original presentation of this paper and draw our attentions to references \cite{LM,MO}.


\end{document}